\documentclass[reqno, 12pt]{amsart}
\def\C{{\mathbb C}}

\def\Q{{\mathbb Q}}
\def\Qp{{\mathbb Q}_p}
\def\Z{{\mathbb Z}}
\def\Zp{{\mathbb Z}_p}

\def\Fp{{{\mathbb F}_p}}
\def\F{{\mathbb F}}
\def\Qp {{{\mathbb Q}_p}}
\def\Fq{{{\mathbb F}_q}}
\def\Fr{{{\mathbb F}_r}}

\newtheorem{lemma}{Lemma}

\newtheorem{prop}{Proposition}

\newtheorem{theorem}{Theorem}

\theoremstyle{definition}
\newtheorem{defn}{Definition}   
\newtheorem{conj}{Conjecture}

\theoremstyle{remark}
\newtheorem{rem}{Remark}        
      
\newtheorem{example}{Example}

\begin{document}
\title[Zeroes of $L$-series in characteristic $p$]
{Zeroes of $L$-series in characteristic $p$}
%\subjclass{Primary 11G09}
%\keywords{Drinfeld modules, $T$-modules, characteristic $p$ $L$-series,
%$L$-zeroes, multi-valued operators}
\author{David Goss}
\thanks{This work is in honor of my mother Barbara Goss Alter.}
\address{Department of Mathematics\\The Ohio State University\\231
W.\ $18^{\text{th}}$ Ave.\\Columbus, Ohio 43210}

\email{goss@math.ohio-state.edu}
\date{February 1, 2006}

\begin{abstract}
In the classical theory of $L$-series, the exact order (of zero) at
a trivial zero is easily computed via the functional equation. In the
characteristic $p$ theory, it has long been known that a functional
equation of classical $s\mapsto 1-s$ type could not exist. In fact,
there exist trivial zeroes whose order of zero is ``too high;''
we call such trivial zeroes ``non-classical.'' This class of 
trivial zeroes was originally studied by Dinesh Thakur \cite{th2}
and quite recently, Javier Diaz-Vargas \cite{dv2}. In the examples 
computed it was found that these non-classical trivial zeroes were
correlated with integers having {\it bounded} sum of $p$-adic
coefficients. In this paper we present a general conjecture along these
lines and explain how this conjecture fits in with previous work
on the zeroes of such characteristic $p$ functions. In particular,  
a solution to this conjecture might entail finding the
``correct'' functional equations in finite characteristic.
\end{abstract}

\maketitle

%****************************************************************************
\section{Introduction}\label{intro}
The debt all mathematicians owe to Euler is obvious and universally known.
Euler's instincts and mathematical taste
have had the most profound effect on all subsequent generations of researchers.
Nowhere is this more evident than in Euler's fantastic contributions to 
number theory  and, in his work on number theory (see, for
instance, \cite{du1}), nothing is more fabulous than
Euler's investigation into what we now call the Riemann zeta function
$\zeta (s)=\zeta_{\Q}(s):=\sum_{n=1}^\infty n^{-s}$.
In fact, Euler was the first to appreciate that $\zeta(s)$ had
a functional equation, (see
\cite{ay1} for a wonderful discussion of Euler's insights).

The fact that there are infinitely many primes goes back to Euclid. Euler
gave an elegant refinement of this result by establishing that 
$$\sum_{p~\rm prime}\frac{1}{p}$$
diverges thereby giving some indication of how well spaced the primes
actually are. Euler's proof of this fact uses the ``Euler product'' associated
to $\zeta (s)$ as well as the divergence of the harmonic
series $\zeta (1)=\sum_{n=1}^\infty \frac{1}{n}$.

Euler was also fascinated by the special values of $\zeta (s)$ at all the
other integers, both positive and negative, and he devoted much energy to their
computation. In fact, it was by studying these values that
Euler discovered the above mentioned functional equation.
A modern number theorist cannot read these discoveries of
Euler without wonder at the sheer beauty and audacity of Euler's methods. 
While the rigorous deduction of the functional equation of
$\zeta (s)$ had to wait until
Riemann and the methods of complex analysis, Euler's argument makes 
thrilling use of both convergent {\it and} divergent series (see
Section \ref{discovery} below). 

In the process of evaluating $\zeta (s)$ where $s$ is an integer, Euler
discovered the ``trivial zeroes'' of $\zeta (s)$ at the negative
even integers. These zeroes play a crucial role in Euler's calculations
as they form the numerator of a rational quantity Euler needs to compute
{\it and} where his methods (and those of everybody else since) fail to
compute the denominator in closed form (see Remark \ref{eulerrem1}
and Equation \ref{euler6}).
Moreover, the functional equation of $\zeta (s)$ then allows one
to then show that these trivial zeroes are simple.
As is well-known in modern number theory, 
the functional equation of $\zeta (s)$ is merely the beginning of a vast
undertaking whose goal is to show that {\it all} arithmetically
interesting Dirichlet series ultimately behave in similar ways.

Euler's work on $\zeta (s)$ has always been an inspiration in our
work on $L$-series in characteristic $p$. We will briefly review the
definition of such characteristic $p$ 
functions in Section \ref{quick} where the base ring $A$ of the
theory of Drinfeld modules plays the role of the integers $\Z$ in
classical theory. In particular, in the prototypical case of $A=\Fr[T]$,
$r=p^m,$ $p$ prime,
(which, like $\Z$, is also Euclidean) one is able to
describe a function $\zeta_A(s):=\sum\limits_{n\in A~\rm monic} n^{-s}$ 
(see Chapter 8 of \cite{go2}) which is a very close
cousin of $\zeta (s)$. Indeed, using the period of the Carlitz module
instead of $2 \pi i$, one could readily establish an analog of Euler's
calculation of $\zeta (2n)$ $n=1,2,\cdots$ for those positive $j$ which
are ``$A$-even'' (i.e. $j$ divisible by $r-1=\#A^\ast$). 
At the negative integers,
$-j$ ($j\geq 0$) one obtains divergent sums of the form
$\sum\limits_{n\in A~\rm monic}n^j$ which, upon regrouping according to the
degree of $n$, become finite. When $j$ is again $A$-even, this sum
is $0$ giving a ``trivial zero.''

Thus, on the surface, the special values of $\zeta_A(s)$ behave very similarly
to those of $\zeta_{\Q}(s)$, and so a first hope would be to follow
Euler and guess at ``the functional equation'' for $\zeta_{\Fr[T]}(s)$. This
fails to work for two basic reasons: 1.\ Obviously the mapping $s\mapsto 
1-s$ is a bijection between even and odd integers; this 
fails for $A$-odd and $A$-even numbers when $r\neq 3$; 2.\ Even with
$r=3$ there are {\it two} distinct analogs of Bernoulli numbers in
the characteristic $p$ theory; this results in considerably more complicated
quotients than in classical theory.

This state of affairs persisted until the mid 1990's when two 
seemingly independent developments occurred. In the first \cite{wa1}
(see also \cite{dv1})
Daqing Wan computed the Newton polygons associated to $\zeta_A(s)$
when $r=p$ (later extended to all $r$ by B.\ Poonen and J.\ Sheats \cite{sh1})
thereby establishing that the absolute value of a zero {\it uniquely}
determines it (including the multiplicity of the zero). Therefore the zeroes of $\zeta_{\Fr[T]}$
lie ``on the line'' $\Fr ((1/T))$ and 
are simple. Wan's calculations were prompted by the observation by the present
author that, in some cases at least, the coefficients of $\zeta_A(s)$ go
to $0$ exponentially.  Such a rate of decay is far faster than is necessitated
to simply establish the basic analyticity properties of such a function
and is implied by having the degrees of certain
``special polynomials'' (see Section \ref{quick}) grow logarithmically.
Such logarithmic growth is now known to be a completely general phenomenon 
\cite{boc1}, \cite{go5}.

In the second development, D.\ Thakur \cite{th2} looked into the possibility
that, for general $A$, the trivial zeroes had higher order multiplicities;
such a phenomenon {\it never} happens classically.
In other words, the construction of trivial zeroes comes equipped with
a ``classical'' lower-bound on the order of zero. However as Thakur
found, there are many instances where 
this lower-bound is {\it not} the exact order; such a trivial zero is
called ``non-classical.''
It is totally remarkable, and very important for us, that Thakur's results
on non-classical trivial zeroes
involve having the sum of the $p$-adic digits of the trivial zero be 
{\it bounded}.
Thakur's computations have recently been extended by Javier Diaz-Vargas
(\cite{dv2}).  These basic results will be recalled in Section \ref{calc}.

The results of Wan, Sheats, etc., are clearly a type of ``Riemann hypothesis''
(see \cite{go3}, \cite{go4}) and one wants to be able to put them into a general conjecture
about the zeroes in {\it full} generality for all motives ($\tau$-sheaves,
etc.) and interpolations at all places of the quotient field $k$ of $A$.
Our first attempt to do so \cite{go3} 
imply ignored the trivial zeroes (as they
are ignored classically in the Riemann hypothesis). As was reported in
\cite{go4},
this conjecture was wrong {\it precisely} because of the impact of higher
order trivial zeroes! More precisely, using the topology on the domain
space of our $L$-series, one is able to use higher order trivial zeroes
to inductively construct $p$-adic integers where the conjecture is false.
The construction produces such $p$-adic integers by building up their
canonical $p$-adic expansion out of the expansions 
of trivial zeroes with higher orders.

It is precisely here that the computations of Thakur and Diaz-Vargas now
fit. Indeed, their computations lead naturally to the conjecture (Conjecture
\ref{general1})
that those integers $j$ for which the trivial zero at $-j$ is non-classical
have {\it bounded} sum of their $p$-adic digits. Presenting this conjecture
is the goal of this work and we discuss the conjecture both at
$\infty$ (i.e., the analog of the complex field) {\it and} at the interpolations of our functions at
finite primes. If this conjecture is
true, then it places a limit on how one can construct counter-examples
to our original conjecture. Indeed, it implies that we need not worry
about non-classical trivial zeroes {\it alone} leading to counter-examples as 
they can have no effect on our construction once the sum of the $p$-adic
digits of the integers being used becomes sufficiently large
(see the discussion in Section \ref{general}; in particular, 
Example \ref{general2}). We view this as positive evidence for the conjecture.

As the reader may see, the trivial zeroes play a special role in both the
classical and characteristic $p$ theory.
But what is the right general conjecture on the zeroes in characteristic
$p$? As of now, we do not
know. However, since the functional equation classically is what allows one to
compute the order of a trivial zero,
it seems to us quite reasonable that a proof of
the above conjecture in our case will generate the correct ideas and
techniques.

It is clear that this work owes a great deal to Euler. 
It should also be clear that it owes a great deal to the 
mathematical taste and insight of Dinesh Thakur and Javier Diaz-Vargas.
It is moreover my pleasure to thank Thakur and J.-P.\ Serre for helpful
comments.

\section{Euler's discovery of the functional equation for $\zeta_\Q(s)$}\label{discovery}
Our treatment here follows that of \cite{ay1}; the reader is referred there for
references and any elided details. Let $\zeta(s)$ be the Riemann zeta
function.

After many years of work, Euler succeeded in computing the values
$\zeta (2n)$, $n=1,2\cdots$ in terms of Bernoulli numbers. Euler
then turned his attention to the values $\zeta (s)$ at negative
integers. Of course, Euler did not have analytic continuation to work
with and relied on his instincts for beauty; nevertheless, he got it
right! Euler begins with the very well known expansion
\begin{equation}\label{euler1}
\frac{1}{1-x}=1+x+x^2+x^3+\cdots+ x^n +\cdots\,.
\end{equation}
Clearly this expansion is only valid when $|x|<1$, but that does not
stop Euler. Upon putting $x=-1$, he deduces
\begin{equation}\label{euler2}
1/2=1-1+1-1+1\cdots\,.
\end{equation}
To the modern eye, this is a nonsensical statement about divergent
series; however following in
Euler's bold steps, we won't let that stop us! Indeed, upon  
applying $x(d/dx)$ to 
Equation \ref{euler1} and plugging in $x=-1$, we obtain 
\begin{equation}\label{euler3}
1/4=1-2+3-4+5\cdots \,.\end{equation}
Applying the process again, Euler finds the ``trivial zero''
\begin{equation}\label{euler4}
0=1-2^2+3^2-\cdots\,,\end{equation}
and so on. Obviously, Euler is not working with the values at the
negative integers of $\zeta (s)$ but rather the function 
\begin{equation}\label{euler5}
\zeta^\ast(s):=(1-2^{1-s})\zeta(s)=\sum_{n=1}^\infty (-1)^{n-1}/n^s \,,
\end{equation}
however
this is of little consequence and the zeta-values Euler obtains are exactly
those rigorously obtained much later by Riemann. (In \cite{ay1}, our
$\zeta^\ast(s)$ is denoted $\phi(s)$.)

Nine years later, in \cite{eu1} (N.B.: while \cite{eu1} was published in
1768, it was written in 1749) Euler notices, at least for small $n\geq
2$, that his calculations imply 
\begin{equation}\label{euler6}
\frac{\zeta^\ast(1-n)}{\zeta^\ast(n)}= \begin{cases}\frac{(-1)^{(n/2)+1}(2^n-1)(n-1)!}{(2^{n-1}-1)\pi^n} & \text{if $n$ is even}\\
0& \text{if $n$ is odd.}\end{cases}\end{equation}
Upon rewriting Equation \ref{euler6} using his gamma function
$\Gamma (s)$ and the cosine, Euler then ``hazards'' to conjecture
\begin{equation}\label{euler7}
\frac{\zeta^\ast (1-s)}{\zeta^\ast (s)}=\frac{-\Gamma (s)(2^s-1)\cos (\pi s/2)}{
(2^{s-1}-1)\pi^s}\,,\end{equation}
which translates easily into the functional equation of $\zeta(s)$!

\begin{rem}\label{eulerrem1}
Note the important role played by the trivial zeroes in Equation
\ref{euler6} in that they
render harmless Euler's inability to calculate explicitly $\zeta^\ast (n)$, 
or $\zeta (n)$, at odd integers $>1$.
\end{rem}

But there is still more! The value $\zeta^\ast (1)$ is precisely the alternating
harmonic series $$1-1/2+1/3-1/4\cdots$$ which Euler knows converges to 
$\log 2$; so his calculations tell him that evaluating the left
hand side of Equation \ref{euler7} at $s=1$ gives $\frac{1}{2\log 2}$. 
Euler then takes the limit on the right hand side and obtains the
same value! To Euler, this is ``strong justification'' for his conjecture 
which Riemann much later proved. (This quote is from Euler's paper
\cite{eu1}, the translation is found at the bottom of page 1083 of
\cite{ay1}.)

\section{A quick introduction to $L$-series in characteristic $p$}\label{quick}
We now briefly go over the basic definitions of characteristic $p$
$L$-series. We will present the general definitions but the reader will
lose very little by always assuming $A=\Fq[T]$ in what follows.

Let $k$ be an arbitrary global function field of transcendency degree $1$
and full field of constants $\Fr$. Let $\infty$ be a fixed place of $k$
of degree $d_\infty$ over $\Fr$
and let $\vert ?\vert_\infty$ be the associated absolute value.
Let $A$ be the Dedekind domain of those functions regular outside $\infty$.
It is easy to see that the unit group of $A$ is the set of non-zero
constants and that one has
\begin{equation}\label{inf1}
h_A=d_\infty\cdot h_k\,,\end{equation}
where $d_\infty$ is the degree of $\infty$ and $h_?$ is the respective
class number.

We let $K$ be the completion of $k$ at $\infty$ and $\F_\infty\simeq
\F_{r^{d_\infty}}\subset K$
be the associated finite field.
We let $\pi\in K$ be a uniformizing element so that every non-zero
element $\alpha$ of $K$ may be written
\begin{equation}\label{inf2}
\alpha=\zeta_\alpha\cdot \pi^{n_\alpha}\cdot\langle \alpha \rangle\,
\end{equation}
where $\zeta_\alpha\in \F_\infty^\ast$, $n_\alpha\in \Z$ and
$\langle \alpha \rangle\in U_1(K)=\{x\in K\mid \vert x\vert_\infty=1\}$
has absolute 
value $1$. The elements $\zeta_\alpha$ and $\langle \alpha\rangle$ 
depend on our choice of $\pi$. The element $\zeta_\alpha$ is called
the ``sign of $\alpha$'' and denoted $\text{sgn}(\alpha)$.

\begin{example}\label{inf3}
When $k=\Fr(T)$ and $A=\Fr[T]$, the simplest choice is $\pi=1/T$ so that 
for $n\in A$ monic of degree $d$, one has
\begin{equation}\label{inf4}
n=\pi^{-d}\langle n \rangle=T^d\langle n\rangle\,, \end{equation}
with $\langle n\rangle=n/T^d\equiv 1\pmod{\pi}$.
\end{example}

In general, an element $\alpha\in K^\ast$ is said to be {\it monic} or
{\it positive} if and only if $\text{sgn}(\alpha)=\zeta_\alpha=1$, which clearly depends on the
choice of $\pi$. Notice that the positive elements clearly form a subgroup
of $K^\ast$.

Let $X$ be the smooth projective curve associated to $k$. For any fractional
ideal $I$ of $A$, we let $\deg_k(I)$ be the degree over $\Fr$ of the divisor
associated to $I$ on the affine curve $X-\infty$. For $\alpha\in k^\ast$,
one sets $\deg_k (\alpha)=\deg_k((\alpha))$ where $(\alpha)$ is
the associated fractional ideal; this clearly gives the correct
degree of a polynomial in $\Fr[T]$.

Let $\C_\infty$ be the completion of a fixed algebraic closure $\bar{K}$
of $K$ equipped with the canonical extension of the normalized absolute
value on $K$.  
\begin{defn}\label{inf5}
Set $S_\infty:=\C_\infty^\ast \times \Zp$\,.
\end{defn}

The space $S_\infty$ plays the role of the complex numbers in our theory
in that it is the domain of ``$n^s$.'' Indeed, let
$s=(x,y)\in S_\infty$ and let $\alpha\in k$ be positive. The element
$u=\langle \alpha \rangle-1$ has absolute value $<1$; thus
$\langle \alpha\rangle^y=(1+u)^y$ is easily defined and computed
via the binomial theorem.
\begin{defn}\label{inf6}
We set
\begin{equation}\label{inf7}
\alpha^s:=x^{\deg_k (\alpha)}\langle \alpha \rangle^y\,.
\end{equation} \end{defn}
\noindent
Clearly $S_\infty$ is a group whose operation is written additively.
Suppose
that $j\in \Z$ and $\alpha^j$ is defined in the usual sense of the
canonical $\Z$-action on the multiplicative group. Let
$\pi_\ast\in \C_\infty^\ast$ be a fixed $d_\infty$-th root of $\pi$. 
Set $s_j:=(\pi_\ast^{-j},j)\in S_\infty$. One checks easily that 
Definition \ref{inf6} gives $\displaystyle \alpha^{s_j}=\alpha^j$. 
When there is no chance of confusion, we denote $s_j$ simply by ``$j$.''

In the basic case $A=\Fr[T]$ one can now proceed to define $L$-series in
complete generality. However, in general there are non-principal ideals.
Fortunately there is a canonical and simple procedure to extend
Definition \ref{inf6} to them as follows. 
Let $\mathcal I$ be the group of fractional ideals of the Dedekind
domain $A$ and let ${\mathcal P}\subseteq \mathcal I$ be the subgroup of
principal ideals. Let ${\mathcal P}^+\subseteq \mathcal P$ be the
subgroup of principal ideals which have positive generators. It is a standard
fact that ${\mathcal I}/{\mathcal P}^+$ is a finite abelian group. The 
association 
\begin{equation}\label{inf8}
{\mathfrak h}\in {\mathcal P}^+\mapsto \langle {\mathfrak h}\rangle:= 
\langle \lambda \rangle\,,\end{equation}
 where $\lambda$ is the unique positive generator of $\mathfrak h$,
is obviously a homomorphism from ${\mathcal P}^+$ to $U_1(K)\subset \C_\infty^\ast$.

Let $U_1(\C_\infty)\subset \C_\infty^\ast$ be the group of $1$-units defined
in the obvious fashion. The binomial theorem again shows that
$U_1(\C_\infty)$ is a $\Zp$-module. However, it is also closed under the
unique operation of taking $p$-th roots; as such $U_1(\C_\infty)$ is
a $\Qp$-vector space.

\begin{lemma}\label{inf9}
The mapping ${\mathcal P}^+\to U_1(\C_\infty)$ given by
$\mathfrak h\mapsto \langle {\mathfrak h}\rangle$ has a unique
extension to $\mathfrak I$ (which we also denote by $\langle ?\rangle$).
\end{lemma}
\begin{proof}
As $U_1(\C_\infty)$ is a $\Qp$-vector space, it is a divisible group; thus the
extension follows by general theory. The uniqueness then follows
by the finitude of ${\mathcal I}/{\mathcal P}^+$. \end{proof}

If $s\in S_\infty$ and $I$ as above, we now set
\begin{equation}\label{inf10}
I^s:=x^{\deg_k (I)}\langle I \rangle^y\,.
\end{equation}
Thus if $\alpha\in k$ is positive one sees that $(\alpha)^s$ agrees with
$\alpha^s$ as in Equation \ref{inf7}.

\subsection{Definition of $L$-series}\label{lseries}

Let $G:=\text{Gal}(k^\text{sep}/k)$ be the absolute Galois group of
$k$ where $k^\text{sep}$ is a fixed separable closure of $k$. Let $\bar{\Q}_p$
be a fixed algebraic closure of $\Qp$ and let $\chi\colon G\to \bar{\Q}_p^\ast$
be a homomorphism of Galois type; i.e., $\chi$ factors through the Galois
group $G_1$ of a {\it finite} abelian extension $k_1$ of $k$. Obviously the
image of $\chi$ consists of roots of unity and viewing these as sitting
in $\C$ (via some injection) allows us to think of $\chi$ as also being
complex valued. In particular, for each place $w$ of $k$ (including
$\infty$) one attaches a local factor as follows: 1. The place $w$ is
ramified for $\chi$; in which case the factor is simply $1$. The place $w$ is
unramified; in which case the factor is $(1-\chi(F_w)t)$ where $F_w$ is
the (arithmetic) Frobenius element at $w$.

Let $R_p\subset \bar{\Q}_p$ be the ring of integers with maximal
ideal $M_p$. We fix an injection of $R_p/M_p$ into $\C_\infty$ and so we
now obtain local factors in $\C_\infty[t]$ for which we will use the same
notation  $(1-\chi(F_w)t)$ etc.

\begin{rem}\label{lseries1}
The reader may be wondering why we simply did not use the obvious
reduction $\bar{\chi}\colon G \to (R_p/M_p)^\ast$ to begin with. The answer is that
there are no non-trivial $p$-power roots of unity in characteristic $p$
and so one is hard pressed to get the local factors right. For instance,
in the case $G_1$ is a $p$-group, the reduced homomorphism $\bar \chi$ is the
trivial character. If however, we would simply use the trivial character
to obtain local factors we would be off at the ramified primes. Thus is is
far better to use the characteristic $0$ factors in the above fashion.
\end{rem}

Let $s\in S_\infty$ and $\chi$ as above. 
\begin{defn}\label{lseries2}
We set
\begin{equation}\label{lseries3}
L(\chi,s):=\prod_{\substack{{v\in \text{Spec}(A)}\\v\,\text{unramified}}} 
(1-\chi(F_v)v^{-s})^{-1}\,.\end{equation}
\end{defn}

As in Section 8.9 of \cite{go2}, 
it is known that $L(\chi,s)$ converges on a ``half-plane''
of $S_\infty$ and can be analytically extended to an ``entire'' function
on $S_\infty$. Thus, one can view $L(\chi,s)$ as a continuous $1$-parameter,
where $y\in\Zp$ is the parameter,
family of entire power series in $x^{-1}$ etc.

While we have only discussed abelian $\chi$ here for simplicity of
exposition, it is clear how to proceed in the non-abelian case.

\subsection{Special polynomials}\label{specialpolys}
Let $j$ now be a non-negative integer with $\chi$, as above and let
$s=(x,y)\in S_\infty$.
\begin{defn}\label{special1}
We set
\begin{equation}\label{special2}
z_L(\chi,x,-j):=L(\chi,\pi_\ast^{j}x,-j)\,.\end{equation}
\end{defn}

It is known that $L(\chi,x,-j)$ is a polynomial in $x^{-1}$ and all such
polynomials are called the {\it special polynomials} of $L(\chi,s)$. 
By unraveling the definition of $z_L(\chi,x,-j)$, one sees that the coefficients
of this polynomial lie in 
\begin{equation}\label{special2.5}
\mathcal{O}:=\mathcal{O}_{\bf V}[\zeta]
\end{equation} where $\zeta$ is a primitive $n$-th root of
unity and $n$ is the order of the reduction $\bar\chi$ of $\chi$ (as a finite character) 
and $\mathcal{O}_{\bf V}\subset \C_\infty$ is the ring of integers in the {\it value field} 
generated by the elements $\{I^{s_1}\}$ (see Subsection 8.2 of \cite{go2}).
As mentioned at the end of \cite{go5},
elementary estimates imply that the degree (in
$x^{-1}$) of $L(\chi,x,-j)$ grows logarithmically in $j$.

\begin{rem}\label{special3}
For $L$-series associated to ``$\tau$-sheaves'' etc., the logarithmic
growth of the special polynomials is due to B\"ockle  \cite{boc1}. 
For arbitrary $L$-series associated to representations of Galois type
(not necessarily abelian) one can use B\"ockle's results and the fact that
the Artin Conjecture is true \cite{go1}
for these functions to deduce logarithmic
growth.\end{rem}

\subsection{Trivial zeroes}\label{trivial}
The classical, characteristic $0$ valued, $L$-series associated to $\chi$ 
also attaches a local factor to the prime $\infty$ if it is unramified
for $\chi$. In the case that $\chi$ is non-principal, one knows that 
this classical $L$-series is entire; i.e., is a polynomial in $u=r^{-s}$.
These infinite factors are missing in the definition of our $L$-series
and thereby equip them with {\it trivial zeroes} as we shall explain here.

Let $\psi$ be a Hayes-module (i.e., a sign normalized rank one Drinfeld
module) \cite{hay1}, Section 7 of \cite{go2} associated to a twisting of $\text{sgn}$.
The module $\psi$ analytically arises from a rank one
lattice generated by an element $\xi\in \C_\infty$; one knows that
$\xi^{r^{d_\infty}-1}\in K^\ast$. The extension 
$K_1:=K(\xi)/K$ is a totally ramified abelian extension with 
Galois group $g_\infty$
isomorphic to $\F_\infty^\ast$ via the action on $\xi$.  This local extension is also obtained
by adjoining to $K$ {\it any} non-trivial division point for $\psi$.

Let $W\subset \bar{\Q}_p$ be the Witt ring of $\F_\infty$ and let
$t\colon g_\infty\to \F_\infty^\ast$ be the homomorphism given by the
action of $g_\infty$ on $\xi$ and let $T_\infty\colon g_\infty\to W^\ast$
be the composition of $t$ and the {\it Teichm\"uller character} of $\F_\infty^\ast$. 

We view $T_\infty$ as being extended in the obvious way 
to a character of the absolute
Galois group $G_\infty$ of $K^{\text{sep}}/K$ where $K^{\text{sep}}\subset
\C_\infty$ is the separable closure.

Let $\chi_\infty$ be the local factor at $\infty$ associated to $\chi$ which
we also view as a character on $G_\infty$.
Assume that for some non-negative $j$ the character $\chi_\infty\cdot T_\infty^j$ is
{\it unramified}. Then, as in Theorem 8.12.5 of \cite{go2}, a double congruence implies that
\begin{equation}\label{trivial1}
z_L(\chi,x,-j)/(1-(\chi_\infty\cdot T_\infty^j)(F_\infty)x^{-d_\infty})\in
\mathcal{O}[x^{-1}]
\,,\end{equation}
where $\mathcal O$ is given in Equation \ref{special2.5}. Thus zeroes
of $(1-(\chi_\infty\cdot T_\infty^j)(F_\infty)x^{-d_\infty})$ clearly give
rise to zeroes of the original $L$-series $L(\chi,s)$.
\begin{defn}\label{trivial2}
The zeroes of $L(\chi,s)$ arising from the factor $(1-(\chi_\infty\cdot
T_\infty^j) (F_\infty)x^{-d_\infty})$ are called {\it trivial zeroes for
$\chi$ at $-j$}. 
\end{defn}

\begin{rem}\label{trivial3}
It is clear how to generalize the above construction of trivial zeroes
to arbitrary 
representations of Galois type. For general $L$-series arising from
Drinfeld modules, $t$-modules, $\tau$-sheaves etc., one proceeds 
cohomologically as in \cite{boc1}.\end{rem}

\begin{rem}\label{trivial4}
If $\chi_\infty\cdot T_\infty^j$ is ramified, then the local factor is
$1$ and so no non-trivial information is deduced. In this case, it is
reasonable to expect that $z_L(\chi,x,-j)$ has no 
zeroes of absolute value $1$. \end{rem}

\begin{defn}\label{trivial5}
Let $t$ be a trivial zero for $L(\chi, s)$ at $-j$; so that 
$\pi_\ast^{j}t$ is a root of 
$$(1-(\chi_\infty\cdot T_\infty^j)(F_\infty)x^{-d_\infty})$$
of order $v_0(t)$. Let $v_1(t)$ be the order of $t$ as a zero
of $z_L(\chi,x,-j)$. By (\ref{trivial1}) we know that $v_0(t)\leq v_1(t)$.
If this inequality is strict, then we say that $t$ is {\it non-classical}.
The set of all non-negative 
$j$ such that $L(\chi,s)$ has a non-classical trivial
zero at $-j$ will be called {\it the non-classical set for $L(\chi,s)$}.
\end{defn}

Let $\bar{\chi}$ be the reduction of $\chi$ considered as a homomorphism
into $\C_\infty^\ast$ via our fixed embedding of $R_p/M_p$ into
$\C_\infty$. Let $\F_\chi$ be the finite field generated by the values
of $\bar\chi$ over the base field $\F_p$; obviously $\F_\chi$ is finite
and we let $q_\chi:=p^{e(\chi)}$ be its order. 

The next proposition is implicit in \cite{th2} and \cite{dv2}
\begin{prop}\label{trivial6}
The non-classical set for $L(\chi,s)$ is closed under multiplication by
$q_\chi$.
\end{prop}
\begin{proof} This follow upon applying the $q_\chi$-th power mapping 
to the coefficients.\end{proof}

\subsection{$v$-adic theory and $v$-adic trivial zeroes}\label{vadic}
Let $v$ be a finite prime of $A$ of degree $d_v$ over $\Fr$. Let $k_v$ be the
local field at $v$ with fixed algebraic closure $\bar{k}_v$, equipped with
the canonical topology, and let $\C_v$ be the associated complete field.
As before let $\bf V$ be the value field and $n$ the order of the reduction
of $\chi$.
Fix an embedding $\sigma\colon {\bf V}[\zeta]\to \C_v$ where
$\zeta$ is a primitive $n$-th root. 

Via $\sigma$, the functions $\{z_L(\chi,x,-j)\}_{j=0}^\infty$ can be thought
of as lying in $\C_v[x^{-d_\infty}]$. Upon setting $x=x_\sigma\in \C_v^\ast$,
they interpolate to a continuous
function $L_{\sigma,v}(\chi,x_\sigma,s_\sigma)$ where
$(x_\sigma,s_\sigma)\in \C_v^\ast \times S_{\sigma, v}$; here
$S_{\sigma,v}=\Zp\times \Z/(r^t-1)=\lim\limits_\leftarrow\Z/p^n(r^t-1)$ 
is the inductive limit over $n$ of
$\Z/p^n(r^t-1)$ and where $r^t-1$ is the number of
roots of unity in the extension of $k_v$ generated by the image
of $\sigma$. 

As in the previous subsection, let $j$ be chosen so that
$\chi_\infty\cdot T_\infty^j$ is unramified at $\infty$. The local
factor $(1-(\chi_\infty\cdot T_\infty^j)(F_\infty)x^{-d_\infty})$ obviously
has roots of unity for its zeroes. Thus, these roots have bounded $v$-adic
absolute value (as of course their absolute value is $1$). Their effect on
the Newton Polygons for $L_{\sigma,v}(\chi,x_\sigma,s_\sigma)$ is thus
very limited and so can essentially be ignored. 

However, the process of interpolating $v$-adically precisely {\it kills}
the Euler factor at $v$ in the following manner. Let $\sigma$ be extended
to the natural action on polynomials via its action on the coefficients.
\begin{prop}\label{vadic1}
One has
\begin{equation}\label{vadic2}
L_{\sigma,v}(\chi,x_\sigma,-j)=\sigma\left((1-\chi(F_v)v^jx_\sigma^{-d_v})z_L(\chi,x_\sigma,-j)\right)\,.\end{equation}\end{prop}
\begin{proof} This follows immediately 
as $\lim\limits_{i\to\infty}v^i=0$ in $\C_v$.\end{proof}

We are thus led to a very interesting class of zeroes for $L_{\sigma,v}
(\chi,x_\sigma,s_\sigma)$.

\begin{defn}\label{vadic3}
The zeroes of $1-\sigma(\chi(F_v)v^j)x_\sigma^{-d_v}$ in $\C_v$ are called the
{\it $v$-adic trivial zeroes of $L_{\sigma,v}(\chi,x_\sigma,s_\sigma)$ at $-j\in S_{\sigma,v}$}.
\end{defn}

The $v$-adic trivial zeroes are remarkably similar to their counterparts
in $S_\infty$. The definition of {\it non-classical $v$-adic trivial
zeroes} is now obvious as is the definition of the {\it non-classical
set for $L_{\sigma,v}(\chi,x_\sigma,s_\sigma)$}.

\begin{rem}\label{vadic4}
Actually, our whole construction of $v$-adic trivial zeroes {\em is}
non-classical; indeed we know 
of no analog of our construction of $v$-adic trivial zeroes in the theory
of $p$-adic $L$-series. However, we will continue to 
use ``non-classical'' to refer
to those trivial zeroes whose order is higher than expected.
\end{rem}

The next result is the obvious analog of Proposition \ref{trivial6}
and has the same proof.

\begin{prop}\label{vadic5}
The non-classical set for $L_{\sigma,v}(\chi,x_\sigma,s_\sigma)$ is closed under 
multiplication by $q_\chi$.\end{prop}

\section{The calculations of Thakur and Diaz-Vargas}\label{calc}
Let $\chi=\chi_0$ be the trivial character with constant value $1$.
\begin{defn}\label{calc1}
We call the function $L(\chi,s)$ the {\it zeta function of $A$} and
denote it $\zeta_A(s)$. The $v$-adic interpolations of $\zeta_A(s)$ 
are denoted $\zeta_{\sigma, v}(x_\sigma,s_\sigma)$.\end{defn}
\noindent
Clearly one has $\zeta_A(s)=\sum_I I^{-s}$ where $I$ runs over the
ideals of $A$. 

The trivial zeroes of $\zeta_A(s)$ then occur at the negative integers
$-j\in S_\infty$ where $j\equiv 0\pmod{r^{d_\infty}-1}$; indeed the
local factor at $\infty$ in this case is simply $1-x^{-d_\infty}$.
In the case $A=\Fr[T]$, one can show that these zeroes are simple;
thus the non-classical set for $\zeta_{\Fr[T]}(s)$ is empty.

\begin{rem}\label{calc1.5}
Let $A=\Fr[T]$ and let $s=(x,y)\in S_\infty$. In \cite{wa1}, \cite{sh1}
it is shown that for fixed $y$, a zero of $\zeta_A(x,y)$, has
multiplicity $1$ and is {\it uniquely}
determined by its absolute value; thus {\it all} zeroes are simple
and must lie in $K$.
Suppose now that $\theta$ is a classical Dirichlet character with
classical (complex) $L$-series, $L(\theta,s)$. Let
$$l(\theta,t):=L(\theta,1/2+it)\,.$$
At the end of \cite{go4}, it is shown how the classical functional
equation combined with the action of complex conjugation imply
that the expansion about $t=0$ of $l(\theta,t)$  is, up to
possible multiplication by a non-trivial constant, a {\it real}
power series.  Of course the Riemann hypothesis is equivalent to
having the zeroes of $l(\theta,t)$ be real; so that classical theory
looks quite similar to what was established for $\zeta_A(s)$.
\end{rem}

Now let $r=p=2$.
\begin{example}\label{calc2}
Let $A:=\F_2[T_1,T_2]/(T_1^2+T_1+T_2^3+T_2+1)$. In this case
the quotient field $k$ has genus $1$.
\end{example}
\begin{example}\label{calc3}
Let $A:=\F_2[T_1, T_2]/(T_1^2+T_1+T_2^5+T_2^3+1)$. Here the quotient
field $k$ has genus $2$.\end{example}
In both cases, one finds that $A$ has class number $1$ which implies
that $d_\infty=1$  also. In both cases, $\zeta_A(s)$ will have
trivial zeroes at all negative integers (as $r-1=1$). 

Let $j$ be an integer and let $l_p(j)$ be the sum of its $p$-adic
digits.
After some hand and machine calculations, Thakur \cite{th2} established the
following result. 
\begin{theorem}\label{calc4} Let $A$ be as in Example \ref{calc2}
or Example \ref{calc3}.  Then the order of vanishing of 
$\zeta_A(s)$ at $s=-j$ is $2$ if $l_2(j)\leq g$ where $g$
is the genus of the quotient field $k$.\end{theorem}
Thus, in particular, the non-classical set for $\zeta_A(s)$ is non-empty.

For very small $j$, Thakur also shows the converse to his result.
Thus, for instance, in the case of Example \ref{calc3}, the trivial
zero at $-7$ is simple. His paper contains other such examples. 

In Theorem 5.4.9 of \cite{th3}, Thakur establishes a partial converse
to Theorem \ref{calc4} in the case of Example \ref{calc2}. More precisely
he shows in this case that the
trivial zero at $s=-j$ is simple if $l_2(j)=2$ or
$j\equiv 0, 3, 5~\text{or}~6\pmod{7}$.

In \cite{dv2}, Javier Diaz-Vargas extends these calculations to more general
$A$ where the degree of $\infty$ is $1$ but where $A$ has non-trivial
class number and so our exponentiation of non-trivial ideals
is used. The same general phenomenon appears to hold.

\section{A general conjecture}\label{general}
Let $w$ be any place of $k$, where $k$ is now completely general,
 and consider the non-classical set $N_w$
for the interpolation of $L(\chi,s)$ at $w$ (so that if $w=\infty$,
this interpolation is $L(\chi,s)$ itself). We know from
Propositions \ref{trivial6} and \ref{vadic5} that $N_w$ is closed under
multiplication by $q_\chi$. Extrapolating vastly from the calculations
presented in the previous section we are led to the following conjecture.

\begin{conj}\label{general1}
The non-classical set $N_w$ consists of elements with {\em bounded}
sum of $p$-adic digits.
\end{conj}

Let $C\geq 0$ be some bound. Then, of course, $\{j\mid l_p(j)\leq C\}$ is
a particularly simple set of positive integers which is
closed under multiplication by any power of $p$.

Suppose that one can find infinitely many $-j$ so that the factor
$1-(\chi_\infty\cdot T_\infty^j)(F_\infty)x^{-d_\infty}$ has many zeroes
of the same absolute value (which will clearly happen if $d_\infty>1$).
Then, in \cite{go4}, 
we showed how to construct elements $\alpha\in \Zp$ so that
the power series arising from $L(\chi,x,\alpha)$ had the strange
property that there were infinitely many slopes (of the associated
Newton Polygon) of length greater
than $1$.
One does this inductively by building
up $\alpha$ via its $p$-adic digits; so that in particular $\alpha$
is built inductively of integers $\{\alpha_i\}$ with $l_p(\alpha_i)$
increasing. 

If Conjecture \ref{general1} is true, then such counter-examples 
{\bf cannot} be created out of non-classical trivial zeroes alone.
This is illustrated by the next example.

\begin{example}\label{general2}
Let $A$ be as in Example \ref{calc2} or \ref{calc3} and suppose
that Conjecture \ref{general1} is true in the sense that 
Thakur's result Theorem \ref{calc4} is both necessary
and sufficient. Then one cannot use
the construction in \cite{go4}
to obtain strange $\alpha$ as above. Indeed, once
our approximations have sufficiently large sum of $p$-adic digits
the effect of the non-classical trivial zeroes is negligible.
In fact, let $\alpha_i$ be an approximation to $\alpha$ with $l_p(\alpha_i)>2$.
Then the trivial zero at $-\alpha_i$ must now  be 
simple and so {\em cannot} contribute
a slope to $\alpha$ of length bigger than $1$.
\end{example}

We view Example \ref{general2} as being some ``justification'' for
our conjecture in the Eulerian sense of Section \ref{discovery}.

Conjecture \ref{general1} is very general but clearly not the final
word. One would like to know the exact structure of the non-classical sets
as well as the exact orders of the associated trivial zeroes. Moreover,
one would like to know how the bounds change as the place $w$
varies, etc. Still Conjecture \ref{general1} is a precise
statement that indicates a much deeper theory of the zeroes.

Finally, in this paper we have worked with representations of Galois
type. It is reasonable to ask for a similar conjecture for arbitrary
motives etc. As of this writing, we do not know how to formulate
such a conjecture. However,
the following simple example indicates some possible structure.
\begin{example}\label{general3}
Let $A$ be as in Example \ref{calc2} or \ref{calc3} so that
$p=2$. Let $\psi$
be the Hayes module associated to $A$ and let $L(\psi,s)$ be
its $L$-series. Then $L(\psi,s)=\zeta_A(s-1)$. Thus $j$ is 
non-classical for $L(\psi,s)$ if and only if $j+1$ is non-classical
for $\zeta_A(s)$. Thus Conjecture \ref{general1} 
implies that $l_p(j+1)$ is bounded. Note that clearly $l_p(j+1)$ can be bounded
while $l_p(j)$ goes to infinity (e.g., $j=2^t-1, t=1,2,\ldots$).
\end{example}

It might be that having $l_p(j+1)$ be  bounded instead of $l_p(j)$
is somehow analogous to having
a functional equation of the form $s\mapsto k-s$ classically for $k\neq 1$.

\end{document}